\newcommand{\xbP}{\Pi}

\newcommand{\xbe}{\in}
\newcommand{\xbf}{\phi}

\newcommand{\xbo}{\omega}

\newcommand{\xbs}{\sigma}
\newcommand{\xbt}{\tau}

\newcommand{\xCB}{A}

\newcommand{\xCK}{\times}

\newcommand{\xCN}{\neg}
\newcommand{\xCQ}{\emptyset}

\newcommand{\xcA}{\forall}

\newcommand{\xcC}{\not\subseteq}

\newcommand{\xcE}{\exists}

\newcommand{\xcS}{\bigcap}

\newcommand{\xcV}{\bigcup}

\newcommand{\xcc}{\subseteq}

\newcommand{\xcl}{\vdash}
\newcommand{\xcm}{\models}

\newcommand{\xcp}{\rightarrow}

\newcommand{\xcr}{\leftrightarrow}
\newcommand{\xcs}{\cap}
\newcommand{\xcu}{\wedge}
\newcommand{\xcv}{\cup}

\newcommand{\xcz}{\Box}

\newcommand{\xda}{{\cal A}}
\newcommand{\xdb}{{\cal B}}

\newcommand{\xdl}{{\cal L}}

\newcommand{\xdx}{{\cal X}}
\newcommand{\xdy}{{\cal Y}}

\newcommand{\xEI}{\begin{itemize}}
\newcommand{\xEJ}{\end{itemize}}

\newcommand{\xEd}{\neq}
\newcommand{\xEh}{\begin{enumerate}}
\newcommand{\xEj}{\end{enumerate}}

\newcommand{\xex}{\lceil}

\newcommand{\Xl}{\ldots}

\newcommand{\bl}{\begin{lemma} \rm}
\newcommand{\el}{\end{lemma}}
\newcommand{\br}{\begin{remark} \rm}
\newcommand{\er}{\end{remark}}
\newcommand{\be}{\begin{example} \rm}
\newcommand{\ee}{\end{example}}
\newcommand{\bco}{\begin{corollary} \rm}
\newcommand{\eco}{\end{corollary}}
\newcommand{\bc}{\begin{claim} \rm}
\newcommand{\ec}{\end{claim}}
\newcommand{\bfa}{\begin{fact} \rm}
\newcommand{\efa}{\end{fact}}
\newcommand{\bp}{\begin{proposition} \rm}
\newcommand{\ep}{\end{proposition}}
\newcommand{\bd}{\begin{definition} \rm}
\newcommand{\ed}{\end{definition}}
\newcommand{\bcs}{\begin{construction} \rm}
\newcommand{\ecs}{\end{construction}}
\newcommand{\bcd}{\begin{condition} \rm}
\newcommand{\ecd}{\end{condition}}
\newcommand{\bt}{\begin{theorem} \rm}
\newcommand{\et}{\end{theorem}}
\newcommand{\bn}{\begin{notation} \rm}
\newcommand{\en}{\end{notation}}
\newcommand{\bfi}{\begin{bild} \rm}
\newcommand{\efi}{\end{bild}}
\newcommand{\bsta}{\begin{statement} \rm}
\newcommand{\esta}{\end{statement}}
\newcommand{\bcom}{\begin{comment} \rm}
\newcommand{\ecom}{\end{comment}}
\newcommand{\bdia}{\begin{diagram} \rm}
\newcommand{\edia}{\end{diagram}}

\newcommand{\bfc}{\begin{figure}[htb] \begin{center}}
\newcommand{\efc}{\end{center} \end{figure}}

\documentstyle[12pt]{article}
\title{FACTORIZATION}

\author{Karl Schlechta
\thanks{
ks@cmi.univ-mrs.fr, karl.schlechta@web.de, http://www.cmi.univ-mrs.fr/ $\sim$ ks
} \\
Laboratoire d'Informatique Fondamentale de Marseille
\thanks{
UMR 6166, CNRS and Universit\'{e} de Provence,
Address: CMI, 39, rue Joliot-Curie, F-13453 Marseille Cedex 13, France
}
}

\date{December 27, 2007}

\begin{document}

\newtheorem{lemma}{Lemma}[section]
\newtheorem{theorem}[lemma]{Theorem}
\newtheorem{proposition}[lemma]{Proposition}
\newtheorem{corollary}[lemma]{Corollary}
\newtheorem{claim}[lemma]{Claim}
\newtheorem{fact}[lemma]{Fact}
\newtheorem{remark}[lemma]{Remark}
\newtheorem{definition}{Definition}[section]
\newtheorem{construction}{Construction}[section]
\newtheorem{condition}{Condition}[section]
\newtheorem{example}{Example}[section]
\newtheorem{notation}{Notation}[section]
\newtheorem{bild}{Figure}[section]
\newtheorem{comment}{Comment}[section]
\newtheorem{statement}{Statement}[section]
\newtheorem{diagram}{Diagram}[section]

\maketitle

\renewcommand{\labelenumi}
  {(\arabic{enumi})}
\renewcommand{\labelenumii}
  {(\arabic{enumi}.\arabic{enumii})}
\renewcommand{\labelenumiii}
  {(\arabic{enumi}.\arabic{enumii}.\arabic{enumiii})}
\renewcommand{\labelenumiv}
  {(\arabic{enumi}.\arabic{enumii}.\arabic{enumiii}.\arabic{enumiv})}

\section{INTRODUCTION}

Parikh and co-authors have investigated a notion of logical independence,
based
on the sharing of essential propositional variables. We do a semantical
analogue here. What Parikh et al. call splitting on the logical level, we
call
factorization (on the semantical level). Note that many of our results are
valid for arbitrary products, not only for classical model sets.

We claim no originality of the basic ideas, just our proofs might be new -
but they are always elementary and very easy.

\paragraph{
The Situation:
}

$\hspace{0.01em}$

We work here with arbitrary, non-empty products. Intuitively, $ \xdy $ is
the set
of models for the propositional variable set $U.$ We assume the Axiom of
Choice.

\bd

$\hspace{0.01em}$


\label{Definition 1.1:}

Let $U$ be an index set, $ \xdy = \xbP \{Y_{k}:k \xbe U\},$ let all $Y_{k}
\xEd \xCQ,$ and $ \xdx \xcc \xdy.$ Thus, $ \xbs \xbe \xdx $
is a function from $U$ to $ \xcV \{Y_{k}:k \xbe U\}$ s.t. $ \xbs (k) \xbe
Y_{k}.$
We then note $X_{k}:=\{y \xbe Y_{k}: \xcE \xbs \xbe \xdx. \xbs (k)=y\}.$

If $U' \xcc U,$ then $ \xbs \xex U' $ will be the restriction of $ \xbs $
to $U',$ and
$ \xdx \xex U':=\{ \xbs \xex U': \xbs \xbe \xdx \}.$

If $ \xda:=\{A_{i}:i \xbe I\}$ is a partition of $U,$ $U' \xcc U,$ then $
\xda \xex U':=\{A_{i} \xcs U' \xEd \xCQ:i \xbe I\}.$

Let $ \xda:=\{A_{i}:i \xbe I\},$ $ \xdb:=\{B_{j}:j \xbe J\}$ both be
partitions of $U,$ then $ \xda $ is called a
refinement of $ \xdb $ iff for all $i \xbe I$ there is $j \xbe J$ s.t.
$A_{i} \xcc B_{j}.$

A partition $ \xda $ of $U$ will be called a factorization of $ \xdx $ iff
$ \xdx =\{ \xbs \xbe \xdy: \xcA i \xbe I( \xbs \xex A_{i} \xbe \xdx \xex
A_{i})\},$ we will also sometimes say for clarity that
$ \xda $ is a partition of $ \xdx $ over $U.$

We will adhere to above notations throughout these pages.

If $ \xdx $ is as above, $U' \xcc U,$ and $ \xbs \xbe \xdx \xex U',$ then
there is obviously some (usually
not unique) $ \xbt \xbe \xdx $ s.t. $ \xbt \xex U' = \xbs.$ This trivial
fact will be used repeatedly in
the following pages. We will denote by $ \xbs^{+}$ some such $ \xbt $ -
context will
tell which are the $U' $ and $U.$ (To be more definite, we may take the
first such $ \xbt $
in some arbitrary enumeration of $ \xdx.)$

Given a propositional language $ \xdl,$ $v( \xdl )$ will be the set of
its
propositional variables, and $v( \xbf )$ the set of variables occuring in
$ \xbf.$
A model set $C$ is called definable iff there is a theory $T$ s.t.
$C=M(T)$ - the
set of models of $T.$

\ed

\section{
THE RESULTS
}

\label{Section THE}

\bfa

$\hspace{0.01em}$


\label{Fact 2.1:}

If $ \xda,$ $ \xdb $ are two partitions of $U,$ $ \xda $ a factorization
of $ \xdx,$
and $ \xda $ a refinement of $ \xdb,$ then $ \xdb $ is also a
factorization of $ \xdx.$

\efa

\paragraph{
Proof:
}

$\hspace{0.01em}$

Trivial by definition. $ \xcz $
\\[3ex]

\bfa

$\hspace{0.01em}$


\label{Fact 2.2:}

Let $ \xda $ be a factorization of $ \xdx $ over $U,$ $U' \xcc U.$ Then $
\xda \xex U' $ is a factorization
of $ \xdx \xex U' $ over $U'.$

\efa

\paragraph{
Proof:
}

$\hspace{0.01em}$

If $A_{i} \xcs U' \xEd \xCQ,$ let $ \xbs '_{i} \xbe \xdx \xex (A_{i} \xcs
U' ).$ Let then $ \xbs_{i}:= \xbs '^{+}_{i} \xex A_{i}.$
If $A_{i} \xcs U' = \xCQ,$ let $ \xbs_{i}:= \xbt \xex A_{i}$ for any $
\xbt \xbe \xdx.$
Then $ \xbs:= \xcV \{ \xbs_{i}:i \xbe I\} \xbe \xdx $ by hypothesis, so $
\xbs \xex U' \xbe \xdx \xex U',$ and $ \xbs \xex (A_{i} \xcs U' \}= \xbs
'_{i}.$
$ \xcz $
\\[3ex]

\bfa

$\hspace{0.01em}$


\label{Fact 2.3:}

If $A \xcv A' $ is a factorization of $ \xdx $ over $U,$ $ \xda $ a
factorization of $ \xdx \xex A$ over $ \xCB,$
$ \xda ' $ a factorization of $ \xdx \xex A' $ over $A',$ then $ \xda
\xcv \xda ' $ is a factorization of $ \xdx $
over $U.$

\efa

\paragraph{
Proof:
}

$\hspace{0.01em}$

Trivial $ \xcz.$
\\[3ex]

\bfa

$\hspace{0.01em}$


\label{Fact 2.4:}

If $ \xda,$ $ \xdb $ are two factorizations of $ \xdx,$ then there is a
common
refining factorization.

\efa

\paragraph{
Proof:
}

$\hspace{0.01em}$

Let $ \xbs $ s.t. $ \xcA i \xbe I \xcA j \xbe J( \xbs \xex (A_{i} \xcs
B_{j}) \xbe \xdx \xex (A_{i} \xcs B_{j})),$ show $ \xbs \xbe \xdx.$
Fix $i \xbe I.$ By Fact 2.2, $ \xdb \xex A_{i}$ is a factorization of $
\xdx \xex A_{i},$ so
$ \xcv \{ \xbs \xex (A_{i} \xcs B_{j}):j \xbe J,$ $A_{i} \xcs B_{j} \xEd
\xCQ \}$ $=$ $ \xbs \xex A_{i} \xbe \xdx \xex A_{i}.$ As $ \xda $ is a
factorization
of $ \xdx,$ $ \xbs \xbe \xdx.$ $ \xcz $
\\[3ex]

This does not generalize to infinitely many factorizations:

\be

$\hspace{0.01em}$


\label{Example 2.1:}

Take as index set $ \xbo +1,$ all $Y_{k}:=\{0,1\}.$
Take $ \xdx:=\{ \xbs:$ $ \xbs \xex \xbo $ arbitrary, and $ \xbs ( \xbo
):=0$ iff $ \xbs \xex \xbo $ is finally constant $\}.$
Consider the partitions $ \xda_{n}:=\{n,( \xbo +1)-n\},$ they are all
fatorizations of $ \xdx,$ as
it suffices to know the sequence from $n+1$ on to know its value on $ \xbo
.$
A common refinement $ \xda $ will have some $A \xbe \xda $ s.t. $ \xbo
\xbe A.$ Suppose there is some
$n \xbe \xbo \xcs A,$ then $A \xcC n+1,$ $A \xcC ( \xbo +1)-(n+1),$ this
is impossible, so $A=\{ \xbo \}.$ If $ \xda $
were a factorization of $ \xdx,$ so would be $\{ \xbo,\{ \xbo \}\}$ by
Fact 2.1, but $ \xdx $
does not factor into $ \xdx \xex \xbo $ and $ \xdx \xex \{ \xbo \}.$

\ee

\bcom

$\hspace{0.01em}$


\label{Comment 2.1:}

Above set $ \xdx $ is not definable as a model set of a corresponding
language $ \xdl:$
If $ \xbf $ is not a tautology, there is a model $m$
s.t. $m \xcm \xCN \xbf.$ $ \xbf $ is finite, let its variables be among
$p_{1}, \Xl,p_{n}$ and perhaps
$p_{ \xbo.}$ If $p_{ \xbo }$ is not among its variables, it is trivially
also false in some $m' $
in $ \xdx.$ If it is, then modify $m$ accordingly beyond $n.$ Thus,
exactly all
tautologies are true in $ \xdx,$ but $ \xdx \xEd \xdy =$ the set of all $
\xdl -$models.

\ecom

We have, however:

\bfa

$\hspace{0.01em}$


\label{Fact 2.5:}

Let $ \xdx = \xcS \{ \xdx_{m}:m \xbe M\}$ and $ \xdx, \xdx_{m} \xcc \xdy
$ for all $m \xbe M.$

Let $ \xda $ be a partition of $U,$ and a factorization of all $
\xdx_{m}.$

Then $ \xda $ is also a factorization of $ \xdx.$

\efa

\paragraph{
Proof:
}

$\hspace{0.01em}$

Let $ \xbs $ s.t. $ \xcA i \xbe I$ $ \xbs \xex A_{i} \xbe \xdx \xex
A_{i}.$

But $ \xdx \xex A_{i}$ $=$ $( \xcS \{ \xdx_{m}:m \xbe M\}) \xex A_{i}$ $
\xcc $ $ \xcS \{ \xdx_{m} \xex A_{i}:m \xbe M\}:$
Let $ \xbt \xbe \xdx \xex A_{i},$ so by $ \xdx = \xcS \{ \xdx_{m}:m \xbe
M\}$
$ \xbt^{+} \xbe \xdx_{m}$ for all $m \xbe M,$ so $ \xbt \xbe \xdx_{m} \xex
A_{i}$ for all $m \xbe M.$

Thus, $ \xcA i \xbe I, \xcA m \xbe M:$ $ \xbs \xex A_{i} \xbe \xdx_{m}
\xex A_{i},$ so $ \xcA m \xbe M. \xbs \xbe \xdx_{m}$ by prerequisite, so
$ \xbs \xbe \xdx.$ $ \xcz $
\\[3ex]

\bfa

$\hspace{0.01em}$


\label{Fact 2.6:}

Let $A \xcv A' $ be a partition of $U,$ and for all $ \xbs \xbe \xdx \xex
A$ and all
$ \xbt:A' \xcp \xcV \{X_{k}:k \xbe A' \}$ with $ \xbt (k) \xbe X_{k}$ $
\xbs \xcv \xbt \xbe \xdx.$ Then

(1) $A \xcv A' $ is a factorization of $ \xdx $ over $U.$

(2) Any partition $ \xda ' =\{A'_{k}:k \xbe I' \}$ of $A' $ is a
factorization of $ \xdx \xex A' $ over $A'.$

(3) If $ \xda $ is a factorization of $ \xdx \xex A$ over $ \xCB,$ and $
\xda ' $ a partition of $A',$
then $ \xda \xcv \xda ' $ is a factorization of $ \xdx.$

\efa

\paragraph{
Proof:
}

$\hspace{0.01em}$

(1) and (2) are trivial, (3) follows from (1), (2), and Fact 2.3. $ \xcz $
\\[3ex]

\bco

$\hspace{0.01em}$


\label{Corollary 2.7:}

Let $U=v( \xdl )$ for some language $ \xdl.$ Let $ \xdx $ be definable,
and $\{ \xda_{m}:m \xbe M\}$ be a
set of factorizations of $ \xdx $ over $U.$ Then $ \xda:= \xcv \{
\xda_{m}:m \xbe M\}$ is also a
factorization of $ \xdx.$

\eco

\paragraph{
Proof:
}

$\hspace{0.01em}$

Let $ \xdx =M(T).$ Consider $ \xbf \xbe T.$ $v( \xbf )$ is finite,
consider $ \xdx \xex v( \xbf ).$
There are only finitely many different ways $v( \xbf )$ is partitioned by
the
$ \xda_{m},$ let them all be among $ \xda_{m_{0}}, \Xl, \xda_{m_{p}}.$
$M( \xbf ) \xex v( \xbf )$ might not be factorized
by all $ \xda_{m_{0}} \xex v( \xbf ), \Xl, \xda_{m_{p}} \xex v( \xbf ),$
but $M(T) \xex v( \xbf )$ is by Fact 2.2.
By Fact 2.4, $ \xda \xex v( \xbf )$ is a factorization of $M(T) \xex v(
\xbf ).$

Consider now $ \xdx_{ \xbf }:=(M(T) \xex v( \xbf )) \xCK \xbP \{(0,1):k
\xbe v( \xdl )-v( \xbf )\}.$

By Fact 2.6, (1) $\{v( \xbf ),v( \xdl )-v( \xbf )\}$ is a factorization of
$ \xdx_{ \xbf }$ over $v( \xdl ).$

By Fact 2.6, (2) $ \xda \xex (v( \xdl )-v( \xbf ))$ is a factorization of
$ \xdx_{ \xbf } \xex (v( \xdl )-v( \xbf ))$
over $v( \xdl )-v( \xbf ).$

By Fact 2.6, (3) $ \xda $ is a factorization of $ \xdx_{ \xbf }$ over $v(
\xdl ).$

$M(T)= \xcS \{(M(T) \xex v( \xbf )) \xCK \xbP \{(0,1):k \xbe v( \xdl )-v(
\xbf )\}$: $ \xbf \xbe T\},$ so by Fact 2.5,
$ \xda $ is a factorization of $M(T).$

$ \xcz $
\\[3ex]

\bcom

$\hspace{0.01em}$


\label{Comment 2.2:}

Obviously, it is unimportant here that we have only 2 truth values, the
proof
would just as well work with any, even an infinite, number of truth
values.
What we really need is the fact that a formula affects only finitely many
propositional variables, and the rest are free.

\ecom

\br

$\hspace{0.01em}$


\label{Remark 2.8:}

The Hamming distance cooperates well with factorization: Let $T \xcl \xbf
,$ and we
want to revise by $ \xCN \xbf.$ Let $M(T)$ factorize into $ \xCB $ and
$A',$ and let $ \xbf $ not
``concern'' $ \xCB.$ Then any $ \xbf -$model can be made to agree on $ \xCB
$ with some T-model,
with an at least as good Hamming distance. (Proof: Take any $ \xbf
-$model,
modify it on $ \xCB $ as you like, it will still be a $ \xbf -$model.)

\er

Unfortunately, the manner of coding can determine if there is a
factorization,
as can be seen by the following example:

\be

$\hspace{0.01em}$


\label{Example 2.2:}

(1)
p= ``$blue$'', q= ``$round$'', q'= ``$blue$ $iff$ $round$''.

Then

$p \xcu q$ = $blue$ $and$ $round$,  $\neg p \xcu \neg q$ =
$\neg blue$ $and$  $\neg round$

$p \xcu q'$ = $blue$ $and$ $round$,  $\neg p \xcu q'$ =
$ \neg blue$ $and$ $\neg round$

\ee

Thus, both code the same (meta-) situation, the first cannot be
factorized,
the second can.

(2)

More generally, we can code e.g. the non-factorising situation
$\{p \xcu q \xcu r, \xCN p \xcu \xCN q \xcu \xCN r\}$ also using $q' =p
\xcr q,$ $r' =p \xcr r,$ and have then the
factorising situation $\{p \xcu q \xcu r, \xCN p \xcu q' \xcu r' \}.$

(3)

The following situation cannot be made factorising:
$\{p \xcu q,$ $p \xcu \xCN q,$ $ \xCN p \xcu \xCN q\}.$ Suppose there were
some such solution.
Then we need some $p' $ and $q',$ and all 4 possibilities
$\{p' \xcu q',$ $p' \xcu \xCN q',$ $ \xCN p' \xcu q',$ $ \xCN p' \xcu
\xCN q' \}.$ If we do not admit impossible
situations (i.e. one of the 4 possibilities is a contradictory coding),
then
2 possibilities have to contain the same situation, e.g. $p \xcu q.$ But
they
are mutually exclusive (as they are negations), so this is impossible.

$ \xcz $
\\[3ex]

\end{document}